\documentclass[11pt]{article}

\usepackage[iso-8859-7]{inputenc}
\usepackage[centertags]{amsmath}
\usepackage{amsfonts}
\usepackage{epsfig}
\usepackage{amssymb}
\usepackage{amsthm}
\usepackage{newlfont}

 \newtheorem{theorem}{Theorem}
 \newtheorem{corollary}[theorem]{Corollary}
 \newtheorem{lemma}[theorem]{Lemma}
 \newtheorem{proposition}[theorem]{Proposition}

\newcommand{\De}{{\Delta}}
\newcommand{\de}{{\delta}}
\newcommand{\Hh}{{\cal H}}
\newcommand{\Ss}{{\cal S}}
\newcommand{\Tt}{{\cal T}}
\newcommand{\Gg}{{\cal G}}
\newcommand{\ww}{\textup{w}}
\newcommand{\ii}{\textup{i}}
\newcommand{\adj}{\textup{adj}}
\newcommand{\diag}{\textup{diag}}
\newcommand{\dist}{\textup{dist}}
\newcommand{\rank}{\textup{rank}}
\newcommand{\pset}{{\cal B}(P,\eps,\ww)}
\newcommand{\la}{{\lambda}}

\newcommand{\si}{{\sigma}}
\newcommand{\Si}{{\Sigma}}
\newcommand{\Real}{\mathbb R}
\newcommand{\Complex}{\mathbb C}

\newcommand{\eps}{\varepsilon}

\begin{document}
%%%%%%%%%%%%%%%%%%%%%%%%%%%%%%%%%%%%%%%%%%%%%%%%%

\bibliographystyle{plain}

\title{ On condition numbers of polynomial eigenvalue problems
          with nonsingular leading coefficients}

\author{Nikolaos Papathanasiou\thanks{Department of Mathematics,
        National Technical University of Athens, Zografou Campus,
        15780 Athens, Greece (nipapath@mail.ntua.gr, N.
        Papathanasiou; ppsarr@math.ntua.gr, P. Psarrakos).}
        ~and Panayiotis Psarrakos\footnotemark[1]}

%%%%%%%%%%%%%%%%%%%%%%%%%%%%%%%%%%%%%%%%%%%%%%%%%%%%%
% Authors and running title to go on top of each page
% \pagestyle{myheadings} \markboth{Nikolaos Papathanasiou \and
% Panayiotis Psarrakos}{On Condition Numbers of Matrix Polynomials}
%%%%%%%%%%%%%%%%%%%%%%%%%%%%%%%%%%%%%%%%%%%%%%%%%%%%%%5

\maketitle

\vspace{-1.10cm}
\begin{center}
   {\small Dedicated to the memory of James H. Wilkinson
   (1919--1986)}
\end{center}

\begin{abstract}
In this paper, we investigate condition numbers of eigenvalue
problems of matrix polynomials with nonsingular leading
coefficients, generalizing classical results of matrix
perturbation theory. We provide a relation between the condition
numbers of eigenvalues and the pseudospectral growth rate. We
obtain that if a simple eigenvalue of a matrix polynomial is
ill-conditioned in some respects, then it is close to be multiple,
and we construct an upper bound for this distance (measured in the
euclidean norm). We also derive a new expression for the condition
number of a simple eigenvalue, which does not involve
eigenvectors. Moreover,   %%%%%%%%%%%% a Bauer-Fike type and
an Elsner-like perturbation bound for matrix polynomials is
presented.
\end{abstract}

{\small


\vspace{0.10cm}


}

% -------------------------------------------------------
\section{Introduction}\label{introduction}

The notions of condition numbers of eigenproblems and eigenvalues
quantify the sensitivity of eigenvalue problems
\cite{BLO,DEMMEL,Kar,KPM,MBO,RUHE,SMITH,Stewart71,STWs,TISS,Wil65,Wil72,Wil84}.
They are widely appreciated tools for investigating the behavior
under perturbations of matrix-based dynamical systems and of
algorithms in numerical linear algebra. An eigenvalue problem is
called \textit{ill-conditioned} (resp., \textit{well-conditioned})
if its condition number is sufficiently large (resp., sufficiently
small).

In 1965, Wilkinson \cite{Wil65} introduced the condition number of
a simple eigenvalue $\la_0$ of a matrix $A$ while discussing the
sensitivity of $\la_0$ in terms of the associated right and left
eigenspaces. Two years later, Smith \cite{SMITH} obtained explicit
expressions for certain condition numbers related to the reduction
of matrix $A$ to its Jordan canonical form. In early 1970's,
Stewart \cite{Stewart71} and Wilkinson \cite{Wil72} used the
condition number of the simple eigenvalue $\la_0$ to construct a
lower bound and an upper bound for the distance from $A$ to the
set of matrices that have $\la_0$ as a multiple eigenvalue,
respectively. Recently, the notion of the condition number of
simple eigenvalues of matrices has been extended to multiple
eigenvalues of matrices \cite{BLO,Kar,KPM} and to eigenvalues of
matrix polynomials \cite{KPM,TISS}.

In this article, we are concerned with conditioning for the
eigenvalue problem of a matrix polynomial $P(\la)$ with a
nonsingular leading coefficient, generalizing known results of
matrix perturbation theory \cite{BLO,ELS,Kar,SMITH,Wil72}.
%%%%%%%%%%%%%%%%%%%%%%%%%%%%%%%%%%%%%%%%%%%%%%%%%%%%%%%%%%%%
% (The case of singular leading coefficients will be discussed in a
% future work.)
%%%%%%%%%%%%%%%%%%%%%%%%%%%%%%%%%%%%%%%%%%%%%%%%%%%%%%%%%%%%%%%
In the next section, we give the definitions and the necessary
background on matrix polynomials. In Section \ref{multiple}, we
investigate the strong connection between the condition numbers of
the eigenvalues of $P(\la)$ and the growth rate of its
pseudospectra. This connection allows us to portrait the
abstraction of the condition numbers of eigenvalues. In Section
\ref{distance}, we examine the relation between the condition
number of a simple eigenvalue $\la_0$ of $P(\la)$ and the distance
from $P(\la)$ to the set of matrix polynomials that have $\la_0$
as a multiple eigenvalue. In particular, we see that if the
condition number of $\la_0$ is sufficiently large, then this
eigenvalue is close to be multiple. In Section \ref{con_num}, we
provide a new expression for the condition number of a simple
eigenvalue $\la_0$ of $P(\la)$, which involves the distances from
$\la_0$ to the rest of the eigenvalues of $P(\la)$. Finally, in
Section \ref{BFEthm}, we present an extension of   %%%%%%%  the Bauer-Fike Theorem \cite{BF,STWs} and
the Elsner Theorem \cite{ELS,STWe,STWs} to matrix polynomials.
Simple numerical examples are also given to illustrate our
results.

% -------------------------------------------------------
\section{Preliminaries on matrix polynomials} \label{definitions}

Consider an $n \times n$ matrix polynomial
\begin{equation} \label{eq:polyP}
   P(\la) \,=\, A_m \la^m + A_{m-1} \la^{m-1} + \cdots + A_1 \la + A_0  ,
\end{equation}
where $\la$ is a complex variable and $A_j \in \Complex^{n \times
n}$ ($j=0,1,\dots,m$) with $\det A_m \ne 0$. The study of matrix
polynomials has a long history, especially with regard to their
spectral analysis, which leads to the solutions of higher order
linear systems of differential equations. The suggested references
on matrix polynomials are \cite{GLR82,LT85,Mk}.

A scalar $\,\la_0 \in \Complex\,$ is called an \textit{eigenvalue}
of $P(\la)$ if the system $\,P(\la_0)x=0\,$ has a nonzero solution
$\,x_0\in\Complex^n $, known as a \textit{right eigenvector} of
$P(\la)$ corresponding to $\la_0$. A nonzero vector $y_0 \in
\Complex^n$ that satisfies $y_0^*P(\la_0)=0$ is called a
\textit{left eigenvector} of $P(\la)$ corresponding to $\la_0$.
The set of all eigenvalues of $P(\la)$ is the \textit{spectrum} of
$P(\la)$, $\,\sigma(P) = \left \{ \la\in\Complex : \det P(\la)=0
\right \} ,\,$ and since $\,\det A_m \ne 0,\,$ it contains no more
than $nm$ distinct (finite) elements. The \textit{algebraic
multiplicity} of an eigenvalue $\,\la_0\in\sigma(P)\,$ is the
multiplicity of $\la_0$ as a zero of the (scalar) polynomial $\det
P(\la)$, and it is always greater than or equal to the
\textit{geometric multiplicity} of $\la_0$, that is, the dimension
of the null space of matrix $P(\la_0)$.

%%%%%%%%%%%%%%%%%%%%%%%%%%%%%%%%%%%%%%%%%%%%%%%%%%%%%%%%%%%
\subsection{Jordan structure and condition number of the eigenproblem}

Let $\la_1 , \la_2,$ $\dots , \la_r \in \si(P)$ be the eigenvalues
of $P(\la)$, where each $\la_i$ appears exactly $k_i$ times if and
only if its geometric multiplicity is $k_i$ ($i = 1 , 2 , \dots ,
r$). Suppose also that for an eigenvalue $\la_i \in \si(P)$, there
exist $\,x_{i,1}, x_{i,2} , \dots ,x_{i,s_i} \in \Complex^{n}\,$
with $\, x_{i,1} \ne 0, \,$ such that
\[
      \sum_{j=1}^{\xi} \frac{1}{ (j-1) ! } \,
        P^{(j-1)}(\la_i)\, x_{i,\xi-j+1} \,=\, 0 \; ; \;\;\;
        \xi = 1,2, \dots , s_i  ,
\]
where the indices denote the derivatives of $P(\la)$ and $s_i$
cannot exceed the algebraic multiplicity of $\la_i$. Then the
vector $x_{i,1}$ is clearly an eigenvector of $\la_i$, and the
vectors $\,x_{i,2} , x_{i,3} , \dots ,x_{i,s_i}\,$ are known as
\textit{generalized eigenvectors}. The set $\,\{ x_{i,1} , x_{i,2}
, \dots ,x_{i,s_i} \}\,$ is called a \textit{Jordan chain of
length} $s_i$ of $P(\la)$ corresponding to the eigenvalue $\la_i$.
Any eigenvalue of $P(\la)$ of geometric multiplicity $k$ has $k$
maximal Jordan chains associated to $k$ linearly independent
eigenvectors, with total number of eigenvectors and generalized
eigenvectors equal to the algebraic multiplicity of this
eigenvalue.

We consider now the $n \times nm$ matrix $\,X = \left [ x_{1,1} \,
\cdots \, x_{1,s_1} \; x_{2,1} \, \cdots \, x_{r,1} \, \cdots \,
x_{r,s_r} \right ]\,$ formed by maximal Jordan chains of $P(\la)$,
and the associated $nm\times nm$ Jordan matrix $\,J = J_1 \oplus
J_2 \oplus \cdots \oplus J_r , \,$ where each $J_i$ is the Jordan
block that corresponds to the Jordan chain $\{ x_{i,1} , x_{i,2} ,
\dots , x_{i,s_i} \}$ of $\la_i$. Then the $nm\times nm$ matrix
$\, Q  = {\small \left [ \begin{array}{cccc}
         X       \\
         X J     \\
         \vdots  \\
         X J^{m-1} \end{array} \right ] }\,$
is invertible \cite{GLR82}, and we can define $\, Y = Q^{-1}
   {\small  \left[ \begin{array}{cccc}
             0       \\
             \vdots  \\
             0       \\
             A_m^{-1}   \end{array} \right ] }$.
The set $(X,J,Y)$ is called a \textit{Jordan triple} of $P(\la)$,
and satisfies $\,P(\la)^{-1} = X(I \la - J)^{-1}Y \,$ for every
scalar $\la \notin \si(P)$ \cite{GLR82}.
%%%%%%%%%%%%%%%%%%%%%%%%%%%%%%%%%%%%%%%%%%%%%%%%%%%%%%%%%%%%%%%%%
Motivated by the latter equality and \cite{CHU}, we define the
\textit{condition number of the eigenproblem}\footnote{Note that
the definition of the condition number $k(P)$ depends on the
choice of the triple $(X,J,Y)$, but for simplicity, the Jordan
triple will not appear explicitly in the notation.} of $P(\la)$ as
$\, k(P) = \left \| X \right \| \, \left \| Y \right \|,\,$ where
$\| \cdot \|$ denotes the \textit{spectral matrix norm}, i.e.,
that norm subordinate to the euclidean vector norm.

%%%%%%%%%%%%%%%%%%%%%%%%%%%%%%%%%%%%%%%%%%%%%%%%%%%%%%%%%%%
\subsection{Companion matrix}

The \textit{(block) companion matrix} of $P(\la)$ is the
$\,nm\times nm\,$ matrix
\[
   C_P \,=\, \left[ \begin{array}{cccc}
         0 &   I    & \cdots & 0      \\
         0 &   0    & \ddots & \vdots \\
    \vdots & \vdots & \ddots & I      \\
-A_m^{-1}A_0 & -A_m^{-1}A_1 &\cdots & -A_m^{-1}A_{m-1}
   \end{array} \right ] .
\]
It is straightforward to verify that
\begin{equation} \label{eq:C_P}
     E(\la)(\la I - C_P)F(\la) \,=\,
         \left[ \begin{array}{cc}
           P(\la) & 0          \\
           0      & I_{m(n-1)}
       \end{array} \right ] ,
\end{equation}
where $\,F(\la) = {\small \left[ \begin{array}{cccc}
         I  & 0      & \cdots &  0       \\
      \la I & I      & \cdots &  0       \\
     \vdots & \vdots & \ddots & \vdots   \\
 \la^{m}I & \la^{m-1}I &  \cdots & I \\
   \end{array} \right ] }$ and $
  \, E(\la) = {\small \left[ \begin{array}{ccccc}
      E_1(\la) & E_2(\la) & \cdots & E_m(\la) \\
        -I & 0 & \cdots & 0          \\
   \vdots & \ddots & \ddots & \vdots \\
         0 &  &  -I & 0              \\
   \end{array} \right ] }$
with $\,E_m(\la) = A_m\,$ and $\,E_r(\la) = A_{r} + \la
E_{r+1}(\la)\,$ for $\,r = m-1, m-2, \dots , 1$. It is also easy
to see that $\det F(\la)=1$ and $\det E(\la)=\pm \det A_m$ ($\ne
0$). As a consequence, $\si(P)$ coincides with the spectrum of
matrix $C_P$, counting algebraic multiplicities.

%%%%%%%%%%%%%%%%%%%%%%%%%%%%%%%%%%%%%%%%%%%%%%%%%%%%%%%%%%%
\subsection{Weighted perturbations and pseudospectrum}

We are interested in perturbations of $P(\la)$ of the form
\begin{equation} \label{eq:polyQ}
  Q(\la) \,=\, P(\la) + \De(\la)
         \,=\, \sum_{j=0}^m (A_j + \De_j) \la^j ,
\end{equation}
where the matrices $\De_0 , \De_1 , \dots , \De_m \in \Complex^{n
\times n}$ are arbitrary. For a given parameter $\eps >0$ and a
given set of nonnegative weights $\ww = \{ w_0, w_1, \dots , w_m
\}$ with $w_0>0$, we define the class of admissible perturbed
matrix polynomials
\[
    \pset = \left \{ Q(\la) \; \mbox{as in} \;(\ref{eq:polyQ}) :
        \| \De_j \| \le \eps\, w_j ,\, j=0,1,\dots ,m \right \}
\]
(recall that $\| \cdot \|$ denotes the spectral matrix norm). The
weights $w_0 , w_1 , \dots , w_m$ allow freedom in how
perturbations are measured, and the set $\pset$ is convex and
compact with respect to the max norm $\,\|P(\la)\|_{\infty}=
\max\limits_{0\le j\le m}\|A_j\|$ \cite{BLP}.

A recently popularized tool for gaining insight into the
sensitivity of eigenvalues to perturbations is pseudospectrum; see
\cite{BLP,ET,LP,TH} and the references therein. The
$\eps$-\textit{pseudospectrum} of $P(\la)$ (introduced in
\cite{TH}) is defined by
\[
   \si_{\eps}(P) \,=\, \left \{ \mu \in \si (Q) :\, Q(\la)
                      \in \pset \right \}
                 \,=\, \left \{ \mu \in \Complex :\,
                 s_{\min}(P(\mu)) \le \eps\,w(|\mu|) \right \} ,
\]
where $s_{\min}(\cdot)$ denotes the minimum singular value of a
matrix and $\,w(\la) = w_m \la^m + w_{m-1} \la^{m-1}$ $+ \cdots +
w_1 \la + w_0 $. The pseudospectrum $\sigma_{\eps}(P)$ is bounded
if and only if $\,\eps\,w_m < s_{\min}(A_m)$ \cite{LP}, and it has
no more connected components than the number of distinct
eigenvalues of $P(\la)$ \cite{BLP}.

%%%%%%%%%%%%%%%%%%%%%%%%%%%%%%%%%%%%%%%%%%%%%%%%%%%%%%%%%%%
\subsection{Condition number of a simple eigenvalue}

Let $\la_0\in \si(P)$ be a simple eigenvalue of $P(\la)$ with
corresponding right eigenvector $x_0 \in \Complex^n$ and left
eigenvector $y_0 \in \Complex^n$ (where both $x_0$ and $y_0$ are
unique up to scalar multiplications). A normwise \textit{condition
number} of the eigenvalue $\la_0$, originally introduced and
studied in \cite{TISS} (in a slightly different form), is defined
by
\begin{eqnarray}
   k(P,\la_0) &=&    \limsup\limits_{\eps \rightarrow 0} \left \{ \frac{\left|\de
                      \la_0\right|}{\eps}:\, \det Q(\la_0 + \de \la_0)  = 0 , \,
                      Q(\la) \in \pset \right \}
                      \label{eq:k(l,P)}  \\
              &=&     \frac{ w ( \left | \la_0 \right | ) \left
                      \| x_0 \right \| \left \| y_0 \right \| }
                      { \left | y_0^* P'(\la_0) x_0 \right | } .
                      \label{eq:k(l,P)2}
\end{eqnarray}

Since $\la_0$ is also a simple eigenvalue of the companion matrix
$C_P$, we define the \textit{condition number of $\la_0$ with
respect to} $C_P$ as
\begin{equation} \label{eq:cond_comp}
       k(C_P,\la_0) \,=\, \frac{\left \| \chi_0 \right \| \left \|
          \psi_0 \right \|}{\left | \psi_0^* \chi_0 \right | }
\end{equation}
(see \cite{RUHE,SMITH,Wil72,Wil84}), where
\begin{equation} \label{eigenvectors}
   \chi_0 \,=\, \left[ \begin{array}{cccc}
         x_0 \\
         \la_0 x_0  \\
         \vdots \\
         \la_0^{m-1} x_0
   \end{array} \right ] \;\;\; \mbox{and} \;\;\;
  \psi_0 \,=\, \left[ \begin{array}{cccc}
          E_1(\la_0)^*y_0 \\
          E_2(\la_0)^*y_0 \\
          \vdots \\
          E_m(\la_0)^*y_0
   \end{array} \right ]
\end{equation}
are associated right and left eigenvectors of $C_P$ for the
eigenvalue $\la_0$, respectively. By straightforward computations,
we can see that $\,\psi_0^* \chi_0 = y_0^* P'(\la_0) x_0$. This
relation and the definitions (\ref{eq:k(l,P)2}) and
(\ref{eq:cond_comp}) yield the following \cite{MCKEY},
\begin{equation}\label{lem: Mackey}
  k(P,\la_0) \,=\, \frac{ w ( \left | \la_0
                \right | )} { \left \| \chi_0 \right \|
                \left \| \psi_0 \right \|} \, k(C_P,\la_0) .
\end{equation}

%%%%%%%%%%%%%%%%%%%%%%%%%%%%%%%%%%%%%%%%%%%%%%%%%%%%%%%%%%%
\subsection{Condition number of a multiple eigenvalue}

Suppose that $\la_0\in \si(P)$ is a multiple eigenvalue of
$P(\la)$, and that $\,p_0\,$ is the maximum length of Jordan
chains corresponding to $\la_0$. Then we can construct a Jordan
triple of $P(\la)$,
\begin{equation} \label{multipletriple}
   (X,J,Y) \,=\, \left ( \left [ x_{1,1} \, \cdots \, x_{1,p_0} \; x_{2,1}
        \,   \cdots \right ] , \, J_1 \oplus J_2 \oplus \cdots \oplus
        J_{\kappa_0} \oplus \tilde{J} , \, \left [
     \begin{array}{c}
          y^*_{1,p_0} \\
          \vdots      \\
          y^*_{1,1}   \\
          y^*_{2,p_0} \\
          \vdots
     \end{array} \right ]\right ) ,
\end{equation}
where $\,J_1, J_2 , \dots , J_{\kappa_0}\,$ are the $p_0 \times
p_0$ Jordan blocks of $\la_0$, and $\tilde J$ contains all the
Jordan blocks of $\la_0$ of order less than $p_0$ and all the
Jordan blocks that correspond to the rest of the eigenvalues of
$P(\la)$. Moreover, $\,x_{1,1},x_{2,1}, \dots , x_{\kappa_0,1}\,$
are right eigenvectors of $P(\la)$ that correspond to $\,J_1, J_2
, \dots , J_{\kappa_0},\,$ and $\,y_{1,1} , y_{2,1} , \dots ,
y_{\kappa_0,1}\,$ are the associated left eigenvectors. Following
the approach of \cite{BLO,Kar,KPM,MBO} on multiple eigenvalues, we
consider the matrices $\,\hat{X} = \left [ x_{1,1} \; x_{2,1}
                   \, \cdots \, x_{\kappa_0,1} \right ]
                   \in \Complex^{n\times \kappa_0}\,$ and
$\,\hat Y = {\small \left [ \begin{array}{c}
        y^*_{1,1}   \\
        y^*_{2,1}   \\
        \vdots      \\
        y^*_{\kappa_0,1} \end{array} \right ] }
        \,\in\, \Complex^{\kappa_0\times n} ,\,$
and define the \textit{condition number of the multiple
eigenvalue} $\la_0$ by
\begin{equation} \label{conditionmultiple}
   \hat{k}(P,\la_0) \,=\, w(| \la_0 |) \, \| \hat X \,\hat Y \| .
\end{equation}
Notice that since the matrices $\hat{X}$ and $\hat{Y}$ are of rank
$\,\kappa_0\le n,\,$ the product $\,\hat{X}\,\hat{Y}\,$ is nonzero
and $\,\hat{k}(P,\la_0) > 0\,$ (keeping in mind that $\,w_0 >0$).
Moreover, if the eigenvalue $\la_0$ is simple, i.e.,
$\,p_0=\kappa_0=1,\,$ then the definitions (\ref{eq:k(l,P)2}) and
(\ref{conditionmultiple}) coincide \cite{KPM}.

% -------------------------------------------------------
\section{Condition numbers of eigenvalues and pseudospectral growth} \label{multiple}

Consider a matrix polynomial $P(\la)$ as in (\ref{eq:polyP}).
Since the leading coefficient of $P(\la)$ is nonsingular, for
sufficiently small $\,\eps ,\,$ the pseudospectrum
$\sigma_{\eps}(P)$ consists of no more than $nm$ bounded connected
components, each one containing a single (possibly multiple)
eigenvalue of $P(\la)$. By the definition (\ref{eq:k(l,P)}) and
the proof of Theorem 5 of \cite{TISS}, it follows that any small
connected component of $\si_{\eps}(P)$ that contains exactly one
simple eigenvalue $\la_0 \in \si(P)$ is approximately a disc
centered at $\la_0$. Recall that the Hausdorff distance between
two sets $\,\Ss,\Tt\subset \Complex\,$ is
\[
    \Hh (\Ss,\Tt) \,=\, \max \left\{  \sup_{s\in \Ss}
    \inf_{t\in \Tt} |s - t| , \, \sup_{t\in \Tt}
    \inf_{s\in \Ss} |s - t| \right \} .
\]

\begin{proposition} \label{prop:simplecircle}
If $\,\la_0\in\si(P)\,$ is a simple eigenvalue of $P(\la)$, then
as $\,\eps \rightarrow 0,\,$ the Hausdorff distance between the
connected component of $\si_{\eps}(P)$ that contains $\la_0$ and
the disc $\left \{ \mu \in\Complex : \, | \mu - \la_0 | \le
k(P,\la_0)\,\eps \right \}$ is $\,o(\eps)$.
\end{proposition}

Next we extend this proposition to multiple eigenvalues of the
matrix polynomial $P(\la)$, generalizing a technique of \cite{Kar}
for matrices (see also \cite{BLO}).

\begin{theorem}  \label{theorem:multiple}
Suppose that $\la_0$ is a multiple eigenvalue of $P(\la)$ and
$p_0$ is the dimension of the maximum Jordan blocks of $\la_0$.
Then as $\,\eps \rightarrow 0,\,$ the Hausdorff distance between
the connected component of pseudospectrum $\si_{\eps}(P)$ that
contains $\la_0$ and the disc $\left \{ \mu \in\Complex :\, | \mu
- \la_0 | \le ( \hat{k}(P,\la_0) \, \eps )^{1/p_0} \right \}\,$ is
$\,o(\eps^{1/p_0})$.
\end{theorem}

\begin{proof}
Consider the Jordan triple $(X,J,Y)$ of $P(\la)$ in
(\ref{multipletriple}) and the condition number $\hat{k}(P,\la_0)$
in (\ref{conditionmultiple}). For sufficiently small $\,\eps >
0,\,$ the pseudospectrum $\si_{\eps}(P)$ has a compact connected
component $\Gg_{\eps}$ such that $\,\Gg_{\eps} \cap \si (P) = \{
\la_0 \}$. In particular, the eigenvalue $\la_0$ lies in the
(nonempty) interior of $\Gg_{\eps}$; see Corollary 3 and Lemma 8
of \cite{BLP}. Let also $\mu$ be a boundary point of $\Gg_{\eps}$.
Then it holds that
\[
    s_{\min}(P(\mu)) \,=\, \eps\,w(|\mu|)
    \;\;\; \mbox{and} \;\;\;
    P(\mu)^{-1} = X ( I \mu - J )^{-1} Y .
\]

Denote now by $N$ the $p_0 \times p_0$ nilpotent matrix having
ones on the super diagonal and zeros elsewhere, and observe that
\[
   N^{p_0} =\, 0 \;\; \mbox{and} \;\;
   ( I\la-N)^{-1} = \left [
    \begin{array}{ccccc}
  \la^{-1} & \la^{-2}  & \dots & \la^{-p_0} \\
  0 & \la^{-1}  & \dots & \la^{-p_0+1} \\
  \vdots & \vdots & \ddots & \vdots \\
  0 & 0 & \dots & \la^{-1}
  \end{array}
  \right ] =
  \la^{-1} \sum\limits_{j=0}^{p_0-1}(\la^{-1}N)^j
  \;\; (\la \ne 0).
\]
As in \cite{BLO,Kar}, we verify that
\begin{eqnarray*}
  \frac{|\mu-\la_0|^{p_0}}{s_{\min} (P(\mu))} &=& |\mu-\la_0|^{p_0}
  \| P^{-1}(\mu) \|
  \,=\, |\mu-\la_0|^{p_0} \| X (I\mu - J)^{-1}Y \| \\
  &=& |\mu-\la_0|^{p_0} \left \| X\, \diag \left \{
  (I\mu - J_1)^{-1} , \dots ,
  (I\mu - J_{\kappa_0})^{-1} ,
  (I\mu - \hat J)^{-1} \right \} Y \right \|   \\
  &=& {\Big \|} (\mu-\la_0)^{p_0-1} X\, \diag \left \{
  \sum\limits_{j=0}^{p_0-1}((\mu-\la_0)^{-1}N)^j ,  \dots
  \right.  \\
  & & \;\;\; \dots ,
  \left. \sum\limits_{j=0}^{p_0-1}((\mu-\la_0)^{-1}N)^j ,
  (\mu-\la_0)(I\mu - \hat J)^{-1}  \right \} Y {\Big \|} .
\end{eqnarray*}
For each one of the first $\kappa_0$ diagonal blocks, we have
\[
   (\mu-\la_0)^{p_0-1} \sum\limits_{j=0}^{p_0-1}((\mu-\la_0)^{-1}N)^j
   \,=\, N^{p_0-1}  + O(\mu - \la_0) .
\]
Thus, it follows
\begin{eqnarray*}
  \frac{|\mu-\la_0|^{p_0}}{s_{\min} (P(\mu))} &=& \left \| X\,
    \diag \left \{ N^{p_0-1} + O(\mu-\la_0) ,
    \dots , N^{p_0-1} + O(\mu-\la_0) ,
     O((\mu-\la_0)^{p_0}) \right \} Y \right \|    \\
    &=& \left \| X\, \diag \left \{ N^{p_0-1} , \dots ,
    N^{p_0-1}  , 0 \right \} Y \right \| + O(|\mu-\la_0|) \\
    &=& {\Big \|}  \left [ 0 \,\cdots\, x_{1,1} \; 0 \,
    \cdots \, x_{2,1} \; 0 \, \cdots \right ] \, \left [
     \begin{array}{c}
          y_{1,p_0}^*   \\
          \vdots        \\
          y^*_{1,1}     \\
          y_{2,p_0}^*   \\
          \vdots
     \end{array} \right ] {\Big \|} + O(|\mu-\la_0|)  ,
\end{eqnarray*}
where the right eigenvectors $\,x_{1,1}, x_{2,1}, \dots,
x_{\kappa_0,1}\,$ and the rows $\,y^*_{1,1}, y^*_{2,1}, \dots,
y^*_{\kappa_0,1}\,$ lie at positions $\,p_0,\,2 p_0,\dots,
\,\kappa_0 p_0,\,$ respectively. As a consequence,
\[
    \frac{|\mu-\la_0|^{p_0}}{s_{\min} (P(\mu))} \,=\, \|\hat X \, \hat Y \|
    + O(|\mu-\la_0|) ,
\]
or
\[
    \frac{|\mu-\la_0|^{p_0}}{\eps \, w(|\mu|)\,\|\hat X \, \hat Y \|}
     \,=\, 1 + O(|\mu-\la_0|) ,
\]
or
\[
    \frac{|\mu-\la_0|}{( \hat{k}(P,\la_0)\, \eps )^{1/p_0}}
     \,=\, 1 + r_{\eps} ,
\]
%%%%%%%%%%%%%%%%%%%%%%%%%%%%%%%%%%%%%%%%%%%%%%%%%%%%%%%%%%%%%%%%%%%%
% or
% \[
%   |\mu-\la_0|^{p_0} \,=\, \eps \, w(|\mu|)\,\|\hat X \, \hat Y \| +
%    \eps \, w(|\mu|)\, O(|\mu-\la_0|) ,
% \]
%%%%%%%%%%%%%%%%%%%%%%%%%%%%%%%%%%%%%%%%%%%%%%%%%%%%%%%%%%%%%%%%%%%%%
where $\,r_{\eps} \in\Real\,$ goes to $\,0\,$ as $\,\eps
\rightarrow 0$. This means that
\[
     |\mu-\la_0| \,=\, ( \hat{k}(P,\la_0)\, \eps )^{1/p_0}
                       + \, o(\eps^{1/p_0}) .
\]
Since $\,\mu\,$ lies on the boundary $\partial\Gg_{\eps}$, it is
easy to see that the Hausdorff distance between $\Gg_{\eps}$ and
the disc $\left \{ \mu \in\Complex :\, | \mu - \la_0 | \le (
\hat{k}(P,\la_0) \, \eps )^{1/p_0} \right \}\,$ is
$\,o(\eps^{1/p_0})$.
\end{proof}

The above two results indicate how the condition number of an
eigenvalue of $P(\la)$ quantifies the sensitivity of this
eigenvalue. Consider, for example, the matrix polynomial
\[
P(\la) \,=\,
  \left [ \begin{array}{ccc}
    (\la-1)^{2}  & \la-1       &  \la-1  \\
      0          & (\la-1)^{2} &    0    \\
      0          & \la^2-1     &  \la^2-1
    \end{array}  \right ]
\]
with $\,\det(P(\la))=(\la-1)^{5}(\la+1)\,$ and $\,\si(P) = \{ 1,\,
-1\}$. The eigenvalue $\la=1$ has algebraic multiplicity $5$ and
geometric multiplicity $3$, and the eigenvalue $\la=-1$ is simple.
A Jordan triple of $P(\la)$ is given by
\[
  X = \left [ \begin{array}{cccccc}
     1 & 0 &  0 & 0 & 0 & 1 \\
     0 & 1 &  1 & 0 & 1 & 0 \\
     0 & 0 & -1 & 1 & 0 & 2
     \end{array} \right ],\,\,
  J = \left [\begin{array}{cccccc}
     1 & 1 & 0 & 0 & 0 &  0 \\
     0 & 1 & 0 & 0 & 0 &  0 \\
     0 & 0 & 1 & 1 & 0 &  0 \\
     0 & 0 & 0 & 1 & 0 &  0 \\
     0 & 0 & 0 & 0 & 1 &  0 \\
     0 & 0 & 0 & 0 & 0 & -1
    \end{array} \right ], \,\,
  Y = \left [ \begin{array}{ccc}
    0 &  0 &  0.25 \\
    1 &  0 & -0.5  \\
    0 &  1 & -0.5  \\
    0 &  1 &   0   \\
   -1 & -1 &   1   \\
    0 &  0 & -0.25
    \end{array} \right ].
\]
The matrices of the eigenvectors that correspond to the maximum
Jordan blocks of eigenvalue $\la=1$ are $\,  \hat X = {\small
\left [ \begin{array}{cc} 1 & 0 \\ 0 & 1 \\ 0 & -1  \end{array}
\right ] }\,$ and $\, \hat Y = {\small \left [ \begin{array}{ccc}
1 & 0 & -0.5 \\ 0 & 1 & 0 \end{array} \right ]}$. Thus, for the
weights $\,w_0=w_1=w_2=1,\,$ we have $\,\hat{k}(P,1) = w(1)\,\|
\hat X \,\hat Y \| = 4.2426$.

%%%%%%%%%%%%%%%%%%%%%%%%%%%%%%%%%%%%%%%%%%%%%%%%%%%%%%%%
\begin{figure}[ht]
\begin{center}
   \vspace*{-1mm} \mbox{ \hspace*{-2mm}
   \epsfig{file=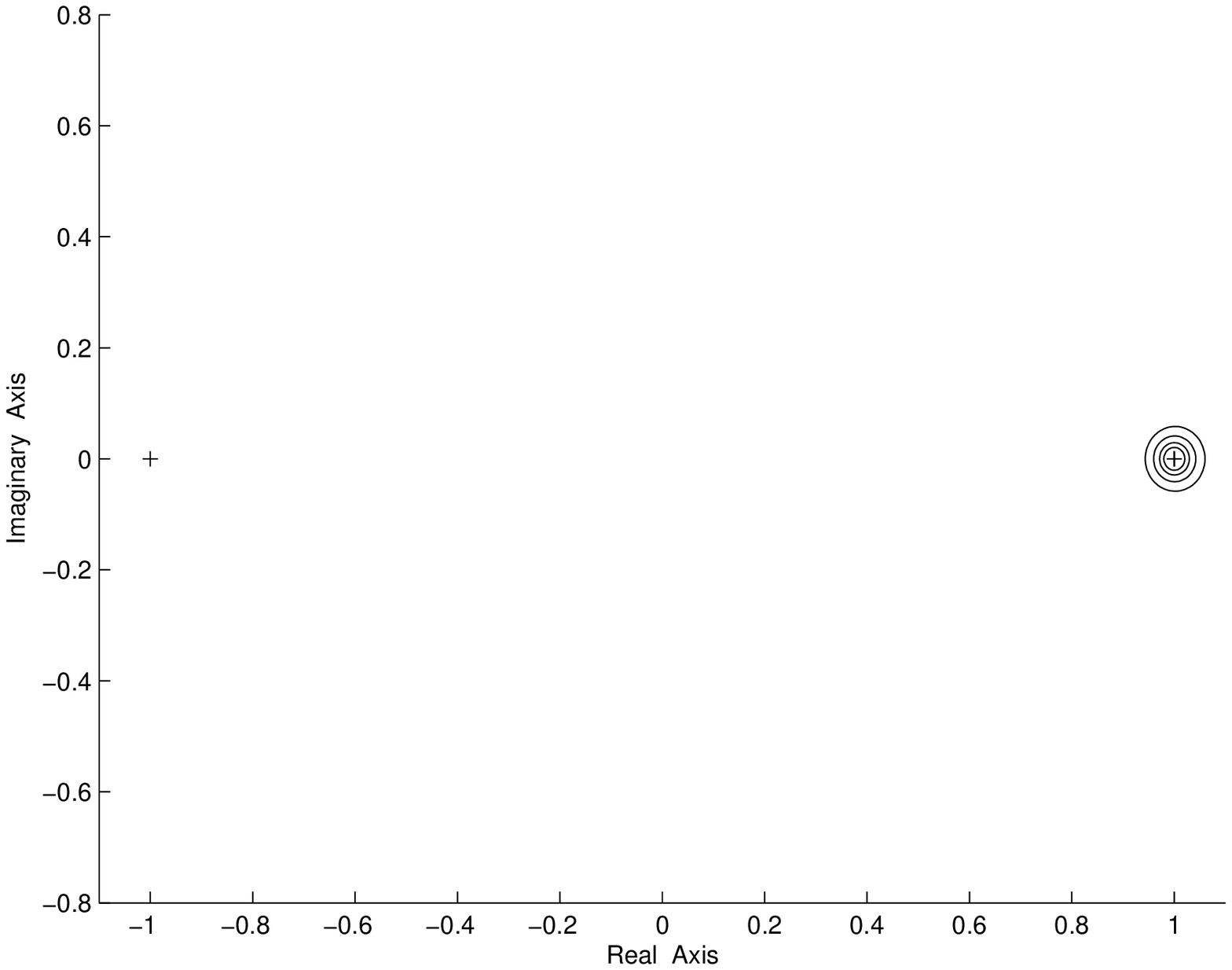, height=56mm, width=68mm, clip=}
   \hspace*{0mm}
   \epsfig{file=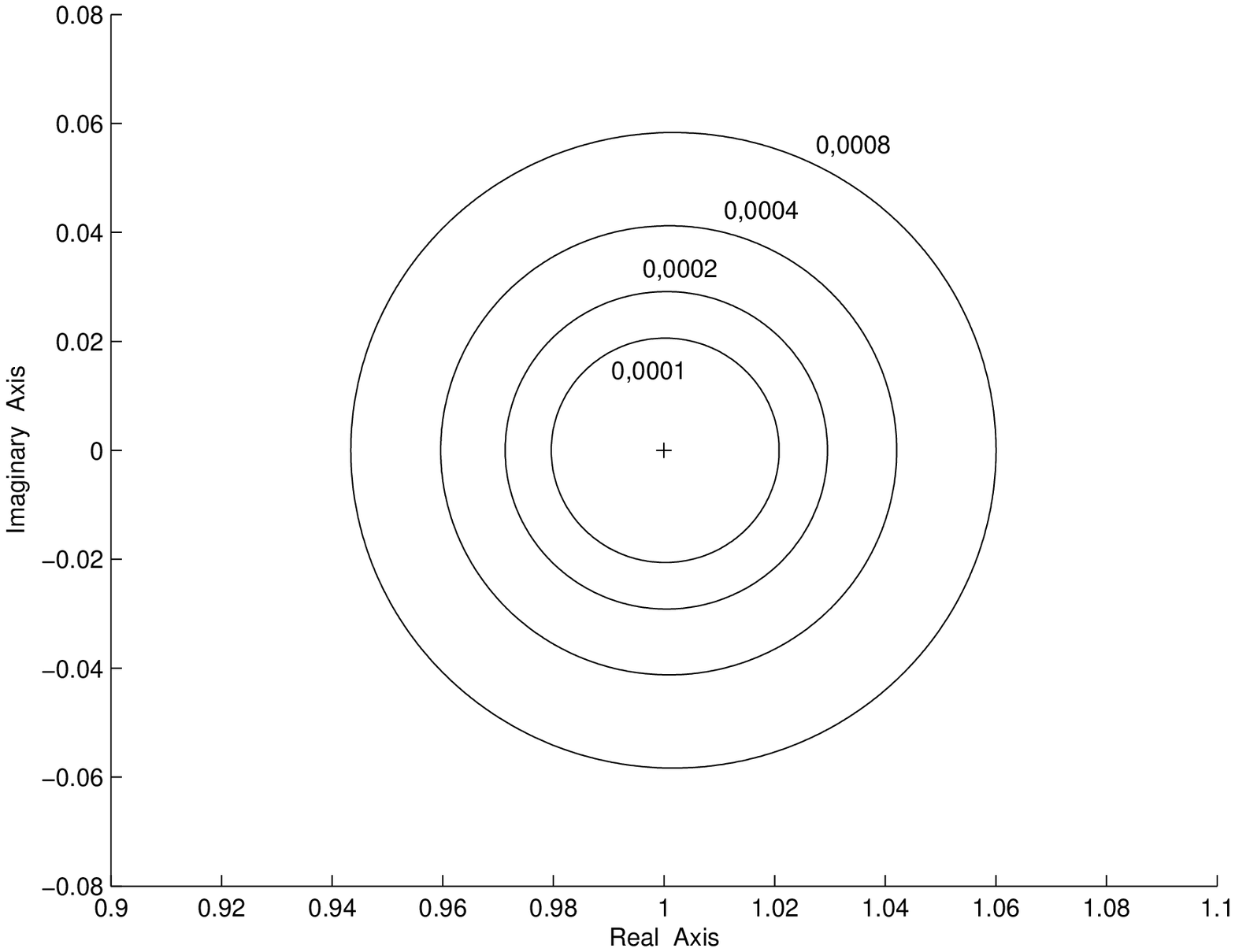, height=56mm, width=68mm, clip=} }
   \vspace{0mm}
   \caption{The boundaries $\,\partial \sigma_{\eps}(P)\,$ for $\,
            \eps = 10^{-4} , 2\cdot 10^{-4} , 4\cdot 10^{-4} , 8\cdot 10^{-4}$.}
   \label{figure2}
\end{center}
\end{figure}
%%%%%%%%%%%%%%%%%%%%%%%%%%%%%%%%%%%%%%%%%%%%%%%%%%%%%%%%
The boundaries of the pseudospectra $\si_{\eps}(P)$, $\,\eps =
10^{-4} , 2\cdot 10^{-4}, 4\cdot 10^{-4}, 8\cdot 10^{-4},\,$ are
illustrated in the left part of Figure \ref{figure2}. The
eigenvalues of $P(\la)$ are marked with $+$'s and the components
of the simple eigenvalue $\la=-1$ are not visible. The components
of the multiple eigenvalue $\la=1$ are magnified in the right part
of the figure, and they are very close to circular discs centered
at $\la = 1$ of radii $\,(\hat{k}(P,1)\, 10^{-4})^{1/2} = 0.0206$,
$(\hat{k}(P,1)\,2\cdot 10^{-4})^{1/2} = 0.0291$,
$(\hat{k}(P,1)\,4\cdot 10^{-4})^{1/2} = 0.0412\,$ and
$\,(\hat{k}(P,1)\,8\cdot 10^{-4})^{1/2} = 0.0583 ,\,$ confirming
Theorem \ref{theorem:multiple}.

% -------------------------------------------------------
\section{Distance from a given simple eigenvalue to multiplicity} \label{distance}

Let $P(\la)$ be a matrix polynomial as in (\ref{eq:polyP}), and
let $\la_0$ be a simple eigenvalue of $P(\la)$. In the sequel, we
generalize a methodology of Wilkinson \cite{Wil72} in order to
obtain a relation between the condition number $k(P,\la_0)$ and
the distance from $P(\la)$ to the matrix polynomials that have
$\la_0$ as a multiple eigenvalue, namely,
\[
   \dist(P,\la_0) \,=\, \inf \left \{ \eps > 0 : \exists \;
      Q(\la) \in \pset \; \mbox{with $\la_0$ as a multiple
      eigenvalue} \right \} .
\]

%%%%%%%%%%%%%%%%%%%%%%%%%%%%%%%%%%%%%%%%%%%%%%%%%%%%%%%%%%%%%%%%
The next proposition is a known result (see \cite[Theorem
3.2]{ANDR} and \cite[Proposition 16]{BLP}). Here, we give a new
proof, which motivates the proof of the main result of this
section (Theorem \ref{thm: Wil_2}) and is necessary for the
remainder.

\begin{proposition} \label{prop: Wil_1}
Let $P(\la)$ be a matrix polynomial as in (\ref{eq:polyP}), $\la_0
\in \si(P)\backslash \si(P')$ and $y_0,x_0 \in \Complex^n$ be
corresponding left and right unit eigenvectors, respectively. If
$\,y_0^* P'(\la_0) x_0$ $= 0 ,\,$ then $\la_0$ is a multiple
eigenvalue of $P(\la)$.
\end{proposition}

\begin{proof}
By Schur's triangularization, and without loss of generality, we
may assume that the matrix $P(\la_0)$ has the following form,
\[
    P(\la_0) = \left[ \begin{array}{cc}
          0 & b^*\\
          0 & B
    \end{array} \right ] \; ; \;\;\;
    b \in \Complex^{n-1}, \,\,
    B \in \Complex^{(n-1) \times (n-1)}  .
\]
Moreover, since $P(\la_0) x_0 = 0$, we can set $x_0 = e_1 = \left
[ \, 1 \;\, 0 \; \cdots \; 0 \, \right ]^T $. Then we have that
$\, y_0^* P'(\la_0) e_1 = 0,\,$ and hence, $\, y_0^* P'(\la_0)
\,=\, [0 \,\, w^*]\,$ for some $\,0 \ne w \in \Complex^{n-1}$.

Since $\la_0 \notin \si(P')$ and $y_0^* P(\la_0) = 0$, it follows
\[
   y_0^* P'(\la_0) \left\{[P'(\la_0)]^{-1}P(\la_0)\right\} \,=\, 0 ,
\]
or equivalently,
\[
y_0^* P'(\la_0) \left\{[P'(\la_0)]^{-1}\left[ \begin{array}{cc}
         0 & b^*\\
         0 & B
   \end{array} \right ]\right\} \,=\, 0 ,
\]
or equivalently,
\[
[0 \,\, w^*]\left[ \begin{array}{cc}
         0 & a^*\\
         0 & A
   \end{array} \right ] \,=\, 0 ,
\]
where $\,a \in \Complex^{n-1}\,$ and $\,A \in \Complex^{(n-1)
\times(n-1)}$. As a consequence, $\,w^* A = 0\,$ and the matrix
$A$ has $0$ as an eigenvalue. Thus, $0$ is a multiple eigenvalue
of the matrix $[P'(\la_0)]^{-1}P(\la_0)$. We consider two cases:

\medskip

\noindent (i) If the geometric multiplicity of $\,0 \in
\si([P'(\la_0)]^{-1}P(\la_0))$ is greater than or equal to $2$,
then $\,\rank(P(\la_0)) \le n-2$, and hence, $\la_0$ is a multiple
eigenvalue of $P(\la)$.

\medskip

\noindent (ii) Suppose that the geometric multiplicity of the
eigenvalue $\,0 \in \si([P'(\la_0)]^{-1} P(\la_0))$ is equal to
$1$ and its algebraic multiplicity is greater than or equal to
$2$. Then, keeping in mind that $\,[P'(\la_0)]^{-1}P(\la_0) e_1 =
0 ,\,$ we verify that there exists a vector $z_1 \in \Complex^n$
such that $\,[P'(\la_0)]^{-1}P(\la_0) z_1 =e_1 ,\,$ or
equivalently, $\, P(\la_0) (-z_1) + P'(\la_0) e_1 = 0$. Thus,
$\la_0$ is a multiple eigenvalue of $P(\la)$ with a Jordan chain
of length at least $2$.
\end{proof}
%%%%%%%%%%%%%%%%%%%%%%%%%%%%%%%%%%%%%%%%%%%%%%%%%%%%%%%%%%%%%%%%%

Recall that the condition number of an invertible matrix $A$ is
defined by $c(A) = \| A \| \, \| A^{-1} \|$ and it is always
greater than or equal to $1$.

\begin{theorem} \label{thm: Wil_2}
Let $P(\la)$ be a matrix polynomial as in (\ref{eq:polyP}), $\la_0
\in \si(P) \backslash \si(P')$ be a simple eigenvalue of $P(\la)$,
and $y_0,x_0 \in \Complex^n$ be corresponding left and right unit
eigenvectors, respectively.
If the vector $[P'(\la_0)]^* y_0$ is not a scalar multiple of $x_0$,
then
\[
  \dist(P,\la_0) \,\le\,
  \frac{ c(P'(\la_0))\, \left\|P(\la_0)\right\|}
  {k(P,\la_0) \left ( \left\|y_0^* P'(\la_0)\right\|^2 -
  | y_0^* P'(\la_0) x_0 |^2 \right )^{1/2}}.
\]
\end{theorem}

\begin{proof}
As in the proof of the previous proposition, without loss of
generality, we may assume that
\[
   P(\la_0) = \left[ \begin{array}{cc}
          0 & b^*\\
          0 & B
    \end{array} \right ] \; ; \;\;\;
    b \in \Complex^{n-1}, \,\,
    B \in \Complex^{(n-1) \times (n-1)}
\]
and $x_0 = e_1$. If we denote $\de = y_0^* P'(\la_0) x_0 = y_0^*
P'(\la_0) e_1 \ne 0$, then it is clear that
\[
    y_0^* P'(\la_0)  \,=\, \left [ \, \de \,\; w^* \right ] ,
\]
for some $\,w \in \Complex^{n-1}$. Furthermore, $\,w \neq 0\,$
because $\,|\de| < \left\|y_0^* P'(\la_0) \right\|$.

Since $\la_0 \notin \si(P')$ and $y_0^* P(\la_0) = 0$, it follows
\[
  y_0^* P'(\la_0) \left \{ [P'(\la_0)]^{-1}P(\la_0) \right \}
  \,=\, 0 ,
\]
or equivalently,
\[
    y_0^* P'(\la_0) \left \{ [ P'(\la_0) ]^{-1} \left [
    \begin{array}{cc}  0 & b^*\\
                       0 & B
   \end{array} \right ] \right \} \,=\, 0 ,
\]
or equivalently,
\[
   \left [ \de \,\, w^* \right ]
   \left[ \begin{array}{cc}
         0 & a^*\\
         0 & A
   \end{array} \right ] \,=\, 0
\]
for some $\,a \in \Complex^{n-1}\,$ and $\,A \in \Complex^{(n-1)
\times(n-1)}$. If $a=0$, then $w^* A = 0$, and the proof of
Proposition \ref{prop: Wil_1} implies that $\la_0$ is a multiple
eigenvalue of $P(\la)$; this is a contradiction. As a consequence,
$a\ne 0$. Moreover,
\[
      w^* A + \de a^* \,=\, 0,
\]
and hence,
\[
   w^* \left(A + \frac{\de}{w^* w} wa^*\right) \,=\, 0 .
\]
This means that if we consider the (perturbation) matrix $\,E =
\left[
\begin{array}{cc} 0 & 0   \\  0 & \frac{\de}{w^* w} wa^*
\end{array} \right ] ,\,$ then the matrix
\[
    [P'(\la_0)]^{-1}P(\la_0) + E \,=\,
    [P'(\la_0)]^{-1} \left [ P(\la_0)
     + P'(\la_0) E \right ]
\]
has $\,0\,$ as a multiple eigenvalue.

We define the $n \times n$ matrices
\[
  \hat \De \,=\, P'(\la_0) E \;\;\; \mbox{and} \;\;\;
  \hat Q = P(\la_0) + \hat \De ,
\]
and the matrix polynomial $\,\De (\la) = \sum_{j=0}^m \De_j
\la^j\,$ with coefficients
\[
 \De_j \,=\, \left ( \frac{ \overline{\la_0} } { |\la_0| } \right )^j
    \frac{w_j}{ w ( | \la_0 | ) } \, \hat \De
    \; ; \;\;\; j = 0 , 1 , \dots , m ,
\]
where (by convention) we assume that $\overline{\la_0} / \la_0  =
0$ whenever $\la_0 = 0$. Then, denoting $\phi = \frac{w'(|\la_0|)}
{w(|\la_0|)} \frac{\overline{\la_0}}{|\la_0|}$, one can verify
that
\[
    \De (\la_0) \,=\, \hat\De
    \;\;\;\mbox{and} \;\;\;
    \De'(\la_0) \,=\, \phi \hat\De .
\]
We define also the matrix polynomial $\,Q(\la) = P(\la) + \De(\la)
,\,$ and consider two cases:

\medskip

\noindent (i) Suppose that the geometric multiplicity of $\,0 \in
\si([P'(\la_0)]^{-1} \hat Q)$ is greater than or equal to $2$.
Then $\,\rank(\hat Q) = \rank(Q(\la_0)) \le n-2,\,$ or
equivalently, $\la_0$ is a multiple eigenvalue of the matrix
polynomial $Q(\la)$ of geometric multiplicity at least $2$.

\medskip

\noindent (ii) Suppose now that the geometric multiplicity of the
eigenvalue $\,0 \in \si([P'(\la_0)]^{-1} \hat Q) \,$ is equal to
$1$, and its algebraic multiplicity is greater than or equal to
$2$. Then, keeping in mind that $\hat Q e_1 = 0$, there is a
vector $z_1 \in \Complex^n$ such that
\[
  [P'(\la_0)]^{-1} \hat Q z_1 \,=\, e_1,
\]
or equivalently,
\begin{equation}\label{eq: Q}
  \hat Q (-z_1) + P'(\la_0) e_1 \,=\, 0.
\end{equation}
We observe that $\,\De'(\la_0) e_1 = \phi \hat\De e_1 = \phi
P'(\la_0) E e_1 = 0$. As a consequence, (\ref{eq: Q}) is written
in the form
\[
    Q(\la_0) (-z_1) + Q'(\la_0) e_1 \,=\, 0 .
\]
Thus, $\la_0$ is a multiple eigenvalue of $Q(\la)$ with a Jordan
chain of length at least $2$.

\medskip

In both cases above, we have proved that $\la_0$ is a multiple
eigenvalue of $Q(\la)$. Furthermore, we see that
\begin{eqnarray*}
      \left \| E \right \|
  &=& \left\|\left[ \begin{array}{cc}
         0 & 0\\
         0 & \frac{\de}{w^* w} wa^*
   \end{array} \right ]\right\| =
   \left\|\frac{\de}{w^* w} wa^*\right\| =
   \frac{|\de|}{\left\|w\right\|}\left\|a\right\| \\
   &\le& \frac{|\de|}{\left\|w\right\|}
   \left\|\left[ \begin{array}{cc}
         0 & a\\
         0 & A
   \end{array} \right ]\right\| =
   \frac{|\de|}{\left\|w\right\|}
   \left\|[P'(\la_0)]^{-1}P(\la_0)\right\|   \\
   &\le& \frac{|\de|}{\left\|w\right\|}
   \left\|[P'(\la_0)]^{-1}\right\|\left\|P(\la_0)\right\| .
\end{eqnarray*}
As a consequence, for every $\,j = 0 , 1 , \dots , m$,
\begin{eqnarray*}
 \left \| \De_j \right \| &=& \frac{w_j}{w(|\la_0|)}\left\|\hat{\De}\right\|
   \,=\,  \frac{w_j}{w(|\la_0|)}\left\|P'(\la_0) E \right \| \\
   &\le& \frac{w_j}{w(|\la_0|)}\left\|P'(\la_0)\right\|
         \left\|E\right\| \\
   &\le& \frac{w_j}{w(|\la_0|)}\frac{|\de|}{ \left \| w \right \|}
         \left \| [P'(\la_0)]^{-1} \right \| \left \| P'(\la_0) \right \|
         \left \| P(\la_0) \right \| \\
%%%%%%%%%%%%%%%%%%%%%%%%%%%%%%%%%%%%%%%%%%%%%%%%%%%%%%%%%%%%%
%   &=& \frac{|\de|}{w(|\la_0|)}\frac{w_j
%       \, c(P'(\la_0)) \left \| P(\la_0) \right \|}
%       {(\left\|y_0^* P'(\la_0)\right\|^2 - \de^2)^{1/2}} \\
%%%%%%%%%%%%%%%%%%%%%%%%%%%%%%%%%%%%%%%%%%%%%%%%%%%%%%%%%%%%%
   &=& w_j \, \frac{ c(P'(\la_0)) \,
       \left \| P(\la_0) \right \| }{k(P,\la_0)
       \left ( \left \| y_0^* P'(\la_0) \right \|^2
       - \de^2 \right )^{1/2} } ,
\end{eqnarray*}
and the proof is complete.
\end{proof}

The spectrum of the matrix polynomial
\begin{eqnarray*}
     P(\la) &=&  I \la^2
                + \left [ \begin{array}{ccc}
                                    1 & 0  & 2    \\
                                    0 & 0  & 0.25 \\
                                    0 & 0  & -0.5
                                  \end{array} \right ] \la
                + \left [ \begin{array}{ccc}
                                    0 & 0  & 8     \\
                                    0 & 25 & -\ii    \\
                                    0 & 0  & 15.25
                                   \end{array} \right ]
\end{eqnarray*}
is $\si(P) = \{0, -1, 0.25 \pm \ii\,3.8971 , \pm \ii\,5\}$. For
the weights $\,w_2 = \left \| A_2\right\| = 1$, $w_1 = \left \|
A_1 \right \| = 2.2919$ and $w_0 = \left \| A_0 \right \| =
25.0379 ,\,$ the above theorem implies $\,\dist(P,-1) \le 0.4991$.
If we estimate the same distance using the method proposed in
\cite{PP}, then we see that $\,\dist(P,-1) \le 0.5991$. On the
other hand, for the eigenvalue $0.25 - \ii\,3.8971$, Theorem
\ref{thm: Wil_2} yields $\,\dist(P,0.25 - \ii\,3.8971) \le
0.1485$, and the method of \cite{PP} implies $\,\dist(P,0.25 -
\ii\,3.8971) \le 0.1398$. At this point, it is necessary to remark
that the methodology of \cite{PP} is applicable to every complex
number and not only to simple eigenvalues of $P(\la)$.

% -------------------------------------------------------

\section{An expression of $k(P,\la_0)$ without eigenvectors} \label{con_num}

In this section, we derive a new expression of the condition
number $k(P,\la_0)$ that involves the distances from $\la_0 \in
\si(P)$ to the rest of the eigenvalues of the matrix polynomial
$P(\la)$, instead of the left and right eigenvectors of $\la_0$.
The next three lemmas are necessary for our discussion. The first
lemma is part of the proof of Theorem 2 in \cite{SMITH}, the
second lemma follows readily from the singular value
decomposition, and the third lemma is part of Theorem 4 in
\cite{SMITH}.

%%%%%%%%%%%%%%%%%%%%%%%%%%%%%%%%%%%%%%%%%%%%%%%%%%%%%%%%%%
% The notation $\adj(\cdot)$ states for the adjugate of a matrix.
%%%%%%%%%%%%%%%%%%%%%%%%%%%%%%%%%%%%%%%%%%%%%%%%%%%%%%%%%%%%%

\begin{lemma} \label{lem: lemSmith15}
For any matrices $\,C,R,W \in \Complex^{n \times n}$, $\, R\, \adj
(W C R) \, W = \det ( W R ) \, \adj (C)$.
\end{lemma}

\begin{lemma}\label{lem: lemSmith2}
Let $A$ be an $n \times n$ matrix with $0$ as a simple eigenvalue,
$\,s_1 \geq s_2 \geq \dots \geq s_{n-1} > s_n = 0 \,$ be the
singular values of $A$, and $u_n , v_n \in \Complex^n$ be left and
right singular vectors of $s_n = 0$, respectively. Then $u_n$ and
$v_n$ are also left and right eigenvectors of $A$ corresponding to
$0$, respectively.
\end{lemma}

\begin{lemma} \label{lem: lemSmith1}
Let $A$ be an $n \times n$ matrix with $0$ as a simple eigenvalue.
If $\,s_1 \geq s_2 \geq \dots \geq s_{n-1} > s_n = 0 \,$ are the
singular values of $A$, then $\,\left \| \adj (A) \right \| = s_1
s_2 \cdots s_{n-1}$.
\end{lemma}

The following theorem is a direct generalization of Theorem 2 of
\cite{SMITH}.

\begin{theorem} \label{thm: thmSmtih}
Let $P(\la)$ be a matrix polynomial as in (\ref{eq:polyP}) with spectrum
$\si(P) = \{\la_1 , \la_2 , \dots, \la_{nm} \}$, counting algebraic
multiplicities. If $\la_i$ is a simple eigenvalue, then
\[
  k(P,\la_i) \,=\,    \frac{w(\left|\la_i\right|)\,
                      \left \| \adj (P(\la_i)) \right\|}
                 {|\det A_m| \, \prod_{j \neq i} |\la_j - \la_i |} .
\]
\end{theorem}

\begin{proof}
For the simple eigenvalue $\la_i \in \si(P)$, consider a singular value
decomposition of matrix $P(\la_i)$,
\[
   P(\la_i) \,=\, U\, \Si \,V^* \,=\, U
   \diag \{ s_1 , \dots , s_{n-1} , 0 \} \, V^* .
\]
Then we have
\[
 \left[ \begin{array}{cc}
         U^* & 0 \\
         0   & I_{n(m-1)}
        \end{array} \right ]
   \left[ \begin{array}{cccc}
        P(\la_i) &  0\\
         0  & I_{n(m-1)}
   \end{array} \right ]
   \left[ \begin{array}{cccc}
         V  &  0\\
         0  & I_{n(m-1)}
   \end{array} \right ] =
   \left[ \begin{array}{cccc}
         \Si &   0    \\
          0  & I_{n(m-1)}
   \end{array} \right ] ,
\]
and Lemma \ref{lem: lemSmith15} implies
\begin{equation} \label{eq:inter}
\left[ \begin{array}{cccc}
          V &  0    \\
          0 &  I_{n(m-1)}
   \end{array} \right ]
  \adj \left ( \left[ \begin{array}{cccc}
        \Si &   0 \\
         0  & I_{n(m-1)}
   \end{array} \right ] \right )
   \left[ \begin{array}{cccc}
       U^*  &    0   \\
         0  & I_{n(m-1)}
   \end{array} \right ]
\end{equation}
\[
   \, = \,   \det(U^*V) \,
        \adj \left ( \left [ \begin{array}{cccc}
        P(\la_0) & 0    \\
           0     & I_{n(m-1)}
   \end{array} \right ] \right ) ,
\]
where $|\det(U^*V)| = 1$.

Let $u_n , v_n \in \Complex^n$ be the last columns of $U$ and $V$,
respectively, i.e., they are left and right singular vectors of
the zero singular value of $P(\la_i)$. Then by Lemma \ref{lem:
lemSmith2}, $y_i=u_n$ and $x_i=v_n$ are left and right unit
eigenvectors of $\,\la_i \in \si(P)$, respectively. Let also
$\psi_i$ and $\chi_i$ be the associated left and right
eigenvectors of $C_P$ for the eigenvalue $\la_i$ given by
(\ref{eigenvectors}). Then by (\ref{lem: Mackey}), \cite[Theorem
2]{SMITH}, Lemma \ref{lem: lemSmith15}, (\ref{eq:C_P}) and
(\ref{eq:inter}) (applied in this specific order), it follows
\begin{eqnarray*}
k(P,\la_i) &=& \frac{w(\left|\la_i\right|)}
{\left\|\chi_i\right\|\left\|\psi_i\right\|} \,k(C_P,\la_i) \\
           &=& \frac{w(\left|\la_i\right|)}
               {\left\|\chi_i\right\|\left\|\psi_i\right\|}
    \frac{\left\|\adj(\la_i I - C_P)\right\|}
    {\prod_{j \neq i} |\la_j - \la_i|}   \\
           &=& \frac{w(\left|\la_i\right|)}
               {\left\|\chi_i\right\|\left\|\psi_i\right\|}
   \frac{\left\|F(\la_i)\,\,\adj(E(\la_i)(\la_i I - C_P)F(\la_i))
   \,E(\la_i)\right\|} { |\det(F(\la_i)E(\la_i))| \,
   \prod_{j\neq i} |\la_j - \la_i|}                   \\
           &=& \frac{w(\left|\la_i\right|)}
               {\left\|\chi_i\right\|\left\|\psi_i\right\|}
  \frac{\left\|F(\la_i)\,\,
    \adj \left ( \left[ \begin{array}{cc}
        P(\la_i) &  0 \\
            0    &  I_{n(m-1)}\\
   \end{array} \right ] \right ) E(\la_i)\right\|}
  { |\det(F(\la_i)E(\la_i))| \,
    \prod_{j \neq i} |\la_j - \la_i| } \\
           &=& \frac{w(\left|\la_i\right|)}
               {\left\|\chi_i\right\|\left\|\psi_i\right\|}
  \frac{\left\|F(\la_i)
    \left[ \begin{array}{cc}
         V &  0 \\
          0 & I_{n(m-1)}
   \end{array} \right ]
  \adj \left ( \left[ \begin{array}{cc}
        \Si &   0 \\
         0  & I_{n(m-1)}
   \end{array} \right ] \right )
   \left[ \begin{array}{cc}
       U^* &    0 \\
         0  & I_{n(m-1)}
   \end{array} \right ]E(\la_i)\right\|}
  { |\det A_m| \, \prod_{j \neq i} |\la_j - \la_i|} .
\end{eqnarray*}
Thus,
\begin{equation}\label{eq: eqSmith1}
   k(P,\la_i) \,=\, \frac{w(\left|\la_i\right|)\, \left\|G\right\|}
               {\left \| \chi_i \right \| \left \| \psi_i \right \| \,
               |\det A_m| \,\prod_{j \neq i} |\la_j - \la_i|} ,
\end{equation}
where
\[
    G \,=\, F(\la_i)
    \left[ \begin{array}{cccc}
         V &  0 \\
          0 & I_{n(m-1)}
   \end{array} \right ]
\adj \left ( \left[ \begin{array}{cccc}
        \Si &   0 \\
         0  & I_{n(m-1)}
   \end{array} \right ] \right )
   \left[ \begin{array}{cccc}
       U^* &    0 \\
         0  & I_{n(m-1)}
   \end{array} \right ] E(\la_i) .
\]

Moreover,
\[
   \adj \left ( \left[ \begin{array}{cccc}
        \Si &   0 \\
         0  & I_{n(m-1)}
   \end{array} \right ] \right ) \,=\,
%%%%%%%%%%%%%%%%%%%%%%%%%%%%%%%%%%%%%%%%%%%
%   \left[ \begin{array}{cccc}
%       0  &        &   &  \\
%          & \ddots &   &  \\
%          &        & s_1 s_2 \ldots s_{n-1}  &  \\
%          &        &   &   0_{n(m-1)}
%   \end{array} \right ] \,=\,
%%%%%%%%%%%%%%%%%%%%%%%%%%%%%%%%%%%%%%%%%%%%%
   \left[ \begin{array}{cccc}
         S  & 0   \\
         0  & 0_{n(m-1)}
   \end{array} \right ],
\]
where $S = s_1 s_2 \cdots s_{n-1} \diag \{ 0 , \dots , 0 , 1 \}$.
As a consequence, the matrix $G$ is written
\begin{eqnarray*}
  G &=& F(\la_i)\,\,
    \left[ \begin{array}{cccc}
         V &  0 \\
         0 & I_{n(m-1)}
   \end{array} \right ]
   \left[ \begin{array}{cccc}
        S  & 0  \\
        0  & 0_{n(m-1)}
   \end{array} \right ]
   \left[ \begin{array}{cccc}
       U^*  &    0    \\
         0  &  I_{n(m-1)}
   \end{array} \right ]\,\,E(\la_i) \\
  &=& F(\la_i)\,\,
  \left[ \begin{array}{cccc}
         VSU^* & 0  \\
          0    & 0_{n(m-1)}
   \end{array} \right ]\,\,E(\la_i) \\
  &=& \left[ \begin{array}{cccc}
             V\,S\,U^* &  0 &  \cdots & 0 \\
       \la_i V\,S\,U^* &  0 &  \cdots & 0  \\
        \vdots & \vdots &  \ddots &  \vdots\\
     \la_i^{m-1} V\,S\,U^* &  0 &  \cdots & 0
   \end{array} \right ]\,\,E(\la_i)          \\
  &=& s_1 s_2 \cdots s_{n-1} \left[ \begin{array}{cccc}
             v_n u_n^*  &  0 &  \cdots & 0   \\
       \la_i v_n u_n^*  &  0 &  \cdots & 0   \\
        \vdots & \vdots &  \ddots &  \vdots  \\
       \la_i^{m-1} v_n u_n^* &  0 &  \cdots & 0
   \end{array} \right ]\,E(\la_i)                 \\
%%%%%%%%%%%%%%%%%%%%%%%%%%%%%%%%%%%%%%%%%%%%%%%%%%%%%%%%%%%%%%%
\end{eqnarray*}
\begin{eqnarray*}
%%%%%%%%%%%%%%%%%%%%%%%%%%%%%%%%%%%%%%%%%%%%%%%%%%%%%%%%%%%%%%%
   &=& s_1 s_2 \cdots s_{n-1} \left[ \begin{array}{cccc}
             x_i y_i^* &  0 &  \cdots & 0 \\
       \la_i x_i y_i^* &  0 &  \cdots & 0  \\
        \vdots & \vdots &  \ddots &  \vdots\\
     \la_i^{m-1} x_i y_i^* &  0   &  \cdots & 0
   \end{array} \right ]\,\,E(\la_i)  \\
   &=& s_1 s_2 \cdots s_{n-1} \left[ \begin{array}{cccc}
             x_i y_i^* E_1(\la_i)&  \cdots &  x_i y_i^* E_m(\la_i)  \\
       \la_i x_i y_i^* E_1(\la_i)&  \cdots &  \la_i x_i y_i^* E_m(\la_i) \\
        \vdots & \ddots &  \vdots  \\
     \la_i^{m-1} x_i y_i^* E_1(\la_i)&  \cdots & \la_i^{m-1} x_i y_i^* E_m(\la_i)
   \end{array} \right ]                 \;\;\;\;\;\;\; \\
   &=& s_1 s_2 \cdots s_{n-1} \left[ \begin{array}{cccc}
             x_i  \\
         \la_i x_i  \\
            \vdots  \\
       \la_i^{m-1} x_i
      \end{array} \right ]
 \left[ \begin{array}{cccc}
            E_1(\la_i)^* y_i \\
            E_2(\la_i)^* y_i \\
                \vdots   \\
            E_m(\la_i) ^* y_i
       \end{array} \right ]^*  \\
   &=& s_1 s_2 \cdots s_{n-1} \,( \chi_i \,\psi_i^*) .
\end{eqnarray*}
Hence, by (\ref{eq: eqSmith1}) and Lemma \ref{lem: lemSmith1}, it follows
\[
  k(P,\la_i) \,=\, \frac{w(\left|\la_i\right|)\, (s_1 s_2 \cdots s_{n-1})
                   \, \left \| \chi_i \,\psi_i^*\right\|}
                   {\left\|\chi_i\right\|\left\|\psi_i\right\|
                   |\det A_m| \, \prod_{j \neq i} |\la_j - \la_i|}
  \,=\, \frac{w(\left|\la_i\right|) \, \left\|\adj (P(\la_i)) \right\|}
        {|\det A_m| \, \prod_{j \neq i} |\la_j - \la_i|}
        \, \frac{ \left \| \chi_i \,\psi_i^* \right \| }
        { \left \| \chi_i \right \| \left \| \psi_i \right \| } .
\]
Since $\, \left \| \chi_i \,\psi_i^* \right \| = \left \| \chi_i \right \|
\left \| \psi_i \right \| ,\,$ the proof is complete.
\end{proof}

The next corollary follows readily.

\begin{corollary}\label{cor: corSmith}
Let $P(\la)$ be a matrix polynomial as in (\ref{eq:polyP}) with
spectrum $\si(P) = \{\la_1 , \la_2 , \dots, \la_{nm} \}$, counting
algebraic multiplicities. If $\la_i$ is a simple eigenvalue of
$P(\la)$ with $\,y_i , x_i \in \Complex^n\,$ associated left and
right unit eigenvectors, respectively, then
\[
  \min \limits_{j \ne i} | \la_j - \la_i | \,\le\, \left (
  \frac{ w( \left | \la_i \right | ) \left \|\adj (P(\la_i)) \right\|}
  { k(P,\la_i) \, |\det A_m| } \right )^{\frac{1}{nm-1}}
  \,=\, \left ( \frac{ \left |  y_i^* P'(\la_i) x_i \right |
  \left \| \adj (P(\la_i)) \right\|}{ |\det A_m| } \right)^{\frac{1}{nm-1} } .
\]
Moreover, if the vector $[P'(\la_0)]^* y_0$ is not a scalar
multiple of $x_0$, then
\[
  \dist(P,\la_i) \,\le\,
  \frac{ c(P'(\la_i))\, \left\|P(\la_i)\right\|\,|\det A_m|  }
  {w(\left|\la_i\right|)\, \left \| \adj (P(\la_i)) \right\|
  \left ( \left\|y_i^* P'(\la_i)\right\|^2 -
  | y_i^* P'(\la_i) x_i |^2 \right )^{1/2}} \,
  \prod_{j \neq i} |\la_j - \la_i | .
\]
\end{corollary}

It is remarkable that for the simple eigenvalue $\,\la_i \in
\si(P)$, Theorem \ref{thm: thmSmtih} and the definition
(\ref{eq:k(l,P)2}) yield
\[
     \left | y_i^* P'(\la_i) x_i \right | \,=\,
     \frac{|\det A_m| \, \prod_{j \ne i} |\la_j - \la_i|}
     {\left\|\adj (P(\la_i)) \right\|} \,\ne\, 0 \;\;\;\;
     \left ( \| x_i \| = \| y_i \| = 1 \right ) .
\]
Thus, Proposition \ref{prop: Wil_1} follows as a corollary of
Theorem \ref{thm: thmSmtih}, and the size of the angle between the
vectors $[P'(\la_i)]^*y_i$ and $x_i$ is partially expressed in
algebraic terms such as determinants and eigenvalues. Note also
that $\la_i$ is relatively close to some other eigenvalues of
$P(\la)$ if and only if $k(P,\la_i)$ is sufficiently greater than
the quantity $\,w(\left|\la_i\right|)\, \left\|\adj (P(\la_i))
\right\| \, |\det A_m|^{-1}$. Furthermore, the condition number
$k(P,\la_i)$ is relatively large (and $\la_i$ is an
ill-conditioned eigenvalue) if and only if the product $\prod_{j
\neq i}|\la_j - \la_i|$ is sufficiently less than
$\,w(\left|\la_i\right|)\, \left\|\adj (P(\la_i)) \right\| \,
|\det A_m|^{-1}$.

To illustrate numerically the latter remark, consider the matrix polynomial
\[
     P(\la) \,=\,  \left [ \begin{array}{cc}
                                    0.001   &  0  \\
                                      0     &  1
                                  \end{array} \right ] \la^2
                  + \left [ \begin{array}{cc}
                                    -0.003   &  0    \\
                                        0    &  -7
                                   \end{array} \right ] \la
                  + \left [ \begin{array}{cc}
                                      0.002  &  0.001  \\
                                         0   &  12
                                   \end{array} \right ]
\]
with (well separated) simple eigenvalues $1$, $2$, $3$ and $4$, and set
$\, w_0 = w_1 = w_2 = 1$. Then it is straightforward to see that
for the eigenvalues $\la=1$ and $\la=2$,
\[
     k(P,1) \cong 3000, \;\,
     |2-1|\,|3-1|\,|4-1| = 6
     \;\;  \mbox{and}  \;\;
     \frac{w(1)\, \left \| \adj (P(1)) \right\|}{|\det A_m|}
     \cong 18000  ,
\]
and
\[
     k(P,2) \cong 7000, \;\,
     |1-2|\,|3-2|\,|4-2| = 2
     \;\;  \mbox{and}  \;\;
     \frac{w(2)\, \left \| \adj (P(2)) \right\|}{|\det A_m|}
     \cong 14000  .
\]
On the other hand, for the eigenvalue $\la=4$, we have
\[
     k(P,4) \cong 21.2897 ,  \;\,
     |1-4|\,|2-4|\,|3-4| = 6
     \;\;  \mbox{and}  \;\;
     \frac{w(4)\, \left \| \adj (P(4)) \right\|}{| \det A_m |}
     \cong 127.738 .
\]

%%%%%%%%%%%%%%%%%%%%%%%%%%%%%%%%%%%%%%%%%%%%%%%%%%%%%%%%
\begin{figure}[ht]
\begin{center}
   \vspace*{-1mm} \mbox{ \hspace*{-2mm}
   \epsfig{file=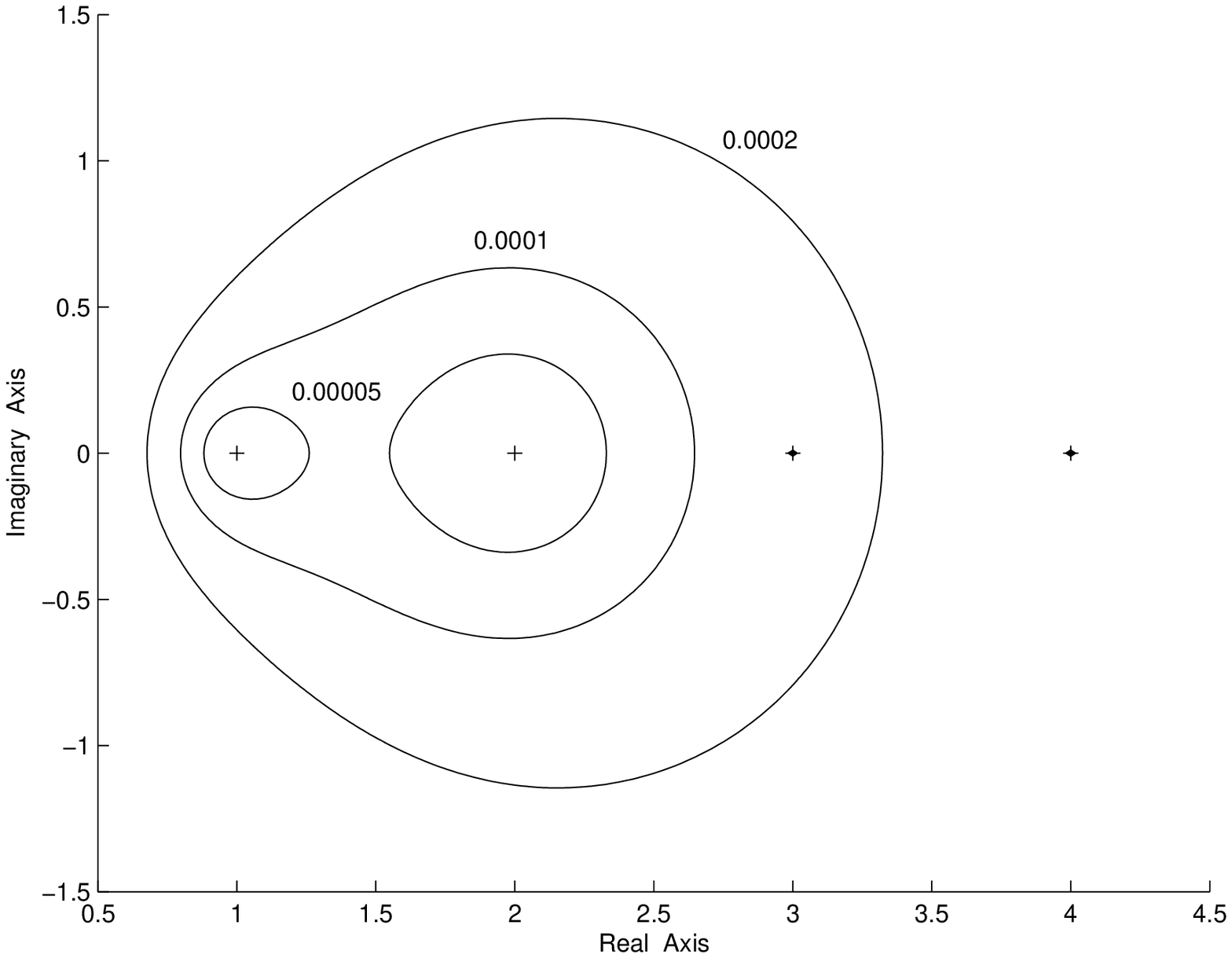, height=56mm, width=68mm, clip=}
   \hspace*{0mm}
   \epsfig{file=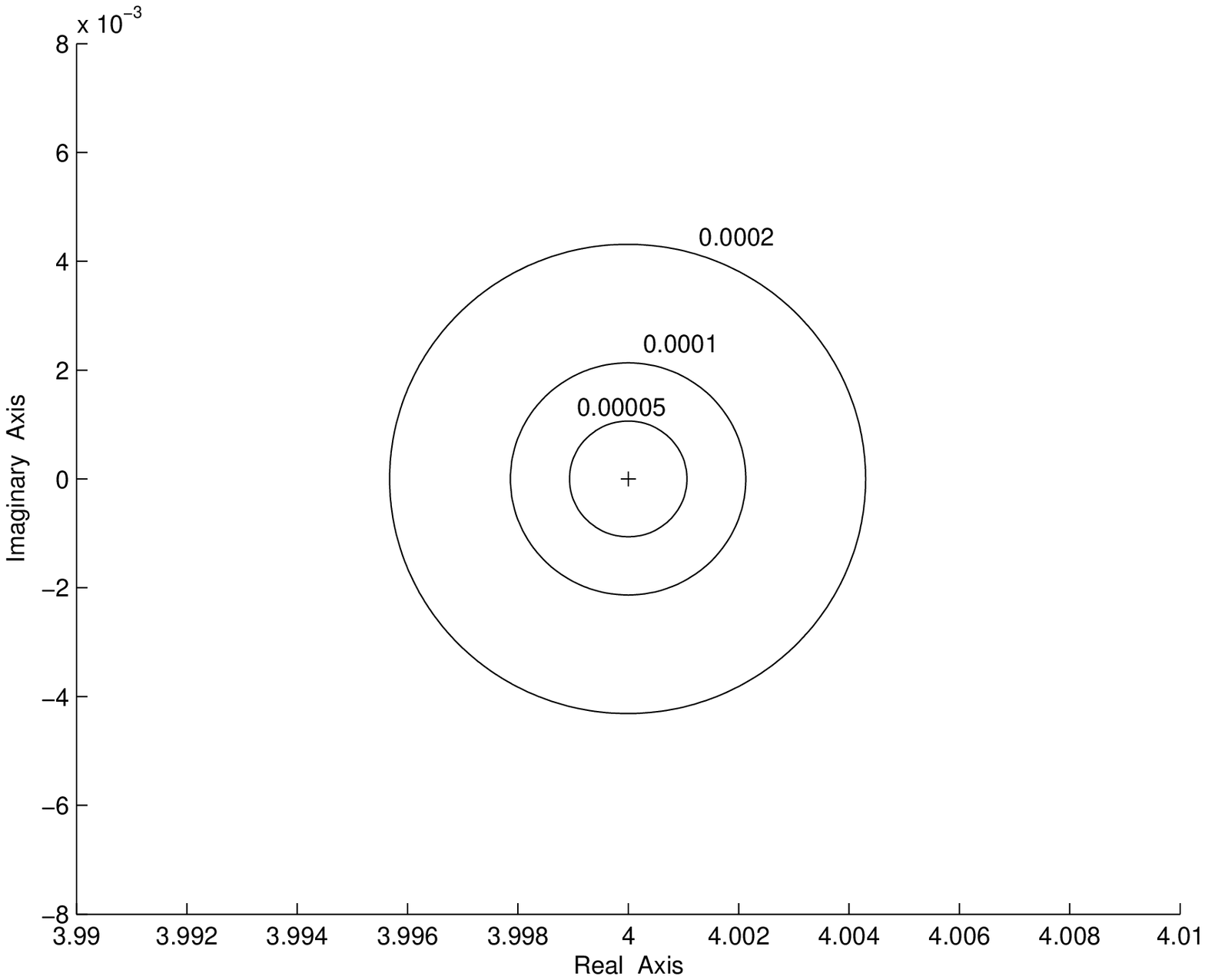, height=56mm, width=68mm, clip=} }
   \vspace{0mm}
   \caption{The boundaries $\,\partial \sigma_{\eps}(P)\,$ for $\,
            \eps = 5\cdot 10^{-5} , 10^{-4} , 2\cdot 10^{-4}$.}
   \label{figure1}
\end{center}
\end{figure}
%%%%%%%%%%%%%%%%%%%%%%%%%%%%%%%%%%%%%%%%%%%%%%%%%%%%%%%%
The left part of Figure \ref{figure1} indicates the boundaries of
the pseudospectra $\si_{\eps}(P)$ for $\,\eps = 5\cdot 10^{-5} ,
10^{-4} , 2\cdot 10^{-4}$. The eigenvalues of $P(\la)$ are marked
with $+$'s. The small components of $\,\si_{5\cdot 10^{-5}} (P)$,
$\si_{10^{-4}} (P)\,$ and $\,\si_{2\cdot 10^{-4}} (P)\,$ that
correspond to the eigenvalue $\la=4$ are not visible in the left
part of the figure, and they are magnified in the right part. Note
that these components almost coincide with circular discs centered
at $\la = 4$ of radii $\,k(P,4)\,\eps$, $\,\eps = 5\cdot 10^{-5} ,
10^{-4} , 2\cdot 10^{-4}$, as expected from Proposition
\ref{prop:simplecircle}. It is also apparently confirmed that the
eigenvalue $\la=2$ is more sensitive than the eigenvalue $\la=1$
(more particularly, one may say that the eigenvalue $\la=2$ is
more than twice as sensitive as $\la=1$), and that both of them
are much more sensitive than the eigenvalue $\la=4$.

% -------------------------------------------------------

\section{An Elsner-like bound} \label{BFEthm}

In this section, we apply the Elsner technique \cite{ELS} (see
also \cite{STWe}) to obtain a perturbation result for matrix
polynomials. This technique allows large perturbations, yielding
error bounds, and it does not distinguish between ill-conditioned
and well-conditioned eigenvalues.

\begin{theorem}\label{thm: thmElsner}
Consider a matrix polynomial $P(\la)$ as in (\ref{eq:polyP}) and a
perturbation $Q(\la) \in \pset$ as in (\ref{eq:polyQ}). For any
$\,\mu \in \si(Q)\backslash\si(P),\,$ it holds that
\[
    \min_{\la\in\si(P)} |\mu - \la| \,\le\,
    \left ( \frac{\eps \, w(|\mu|)}{|\det A_m|}
    \right )^{\frac{1}{mn}}
    \left \|P(\mu) \right\|^{1-\frac{1}{mn}}.
\]
\end{theorem}

\begin{proof} Let $\,\si(P) = \{ \la_1 , \la_2 , \dots , \la_{nm}
\},\,$ counting algebraic multiplicities, and suppose $\mu \in
\si(Q)\backslash\si(P)$. Then
\[
    \min_{\la\in\si(P)} \left | \mu - \la \right |^{nm}
   \,\le\, \prod_{i=1}^{nm} |\mu - \la_i|
   \,=\,   \frac{|\det P(\mu)|}{|\det A_m|} .
\]
Let now $U = \left [ \, u_1 \;\, u_2 \; \cdots \; u_n \, \right ]$
be an $n\times n$ unitary matrix such that $Q(\mu)u_1=0$. By
Hadamard's inequality \cite{STWe} (see also \cite[Theorem
2.4]{STWs}), it follows
\begin{eqnarray*}
    \min_{\la\in\si(P)}\left |\mu - \la \right |^{nm} &\le&
    \frac{|\det P(\mu)|}{|\det A_m|}
    \,=\, \frac{|\det P(\mu)|\,\,|\det U|}{|\det A_m|} \\
    &\le& \frac{1}{|\det A_m|}\,\,\prod_{i=1}^{nm}
          \left \| P(\mu) u_i \right\|       \\
%%%%%%%%%%%%%%%%%%%%%%%%%%%%%%%%%%%%%%%%%%%%%%%%%%%%%%%%%%%%%%%
% \end{eqnarray*}
% \begin{eqnarray*}
%%%%%%%%%%%%%%%%%%%%%%%%%%%%%%%%%%%%%%%%%%%%%%%%%%%%%%%%%%%%%%%
    &=& \frac{1}{|\det A_m|} \left \| P(\mu) u_1 \right \|
        \prod_{i=2}^{nm} \left \| P(\mu) u_i \right \| \\
    &=& \frac{1}{|\det A_m|}\left\|P(\mu)u_1 - Q(\mu)u_1\right\|\,\,
        \prod_{i=2}^{nm} \left \| P(\mu) u_i \right \| \\
    &\le& \frac{1}{|\det A_m|}\left\|\De(\mu)u_1\right\|
    \left\|P(\mu)\right\|^{nm-1} \\
    &\le& \frac{\eps \,\, w(|\mu|)}{|\det A_m|}
    \left\|P(\mu)\right\|^{nm-1} ,
\end{eqnarray*}
and the proof is complete.
\end{proof}

Recently, the classical Bauer-Fike Theorem \cite{BF} has been
generalized to the case of matrix polynomials \cite{CHU}. Applying
the arguments of the proof of Theorem 4.1 in \cite{CHU}, it is
easy to verify the ``weighted version'' of the result.

\begin{theorem}\label{thm: thmChu}
Consider a matrix polynomial $P(\la)$ as in (\ref{eq:polyP}) and a
perturbation $Q(\la) \in \pset$ as in (\ref{eq:polyQ}), and let
$(X, J, Y)$ be a Jordan triple of $P(\la)$. For any $\mu \in
\si(Q)\backslash\si(P)$, it holds that $\,
\min\limits_{\la\in\si(P)} |\mu - \la| \le \max \left\{\vartheta ,
\, \vartheta^{1/p}\right\},\,$ where $\,\vartheta = p \ k(P) \,
\eps \, w(|\mu|)\,$ and $\,p\,$ is the maximum dimension of the
Jordan blocks of $J$.
\end{theorem}

To compare these two bounds, we consider the matrix polynomial
\[
     P(\la) \,=\,  I \la^3
             + \left [ \begin{array}{cc}
                                    0 & \sqrt{2}  \\
                             \sqrt{2} & 0
                                  \end{array} \right ] \la^2
               + \left [ \begin{array}{cc}
                                    0 & -1  \\
                                    1 &  0
                                   \end{array} \right ] \la
\]
(see \cite[Example 1]{CHU} and \cite[Example 1.5]{GLR82})
with $\, \det P(\la) = \la^2 (\la+1)^2 (\la-1)^2 $.
A Jordan triple $(X,J,Y)$ of $P(\la)$ is given by
\[
    X  = \left [\begin{array}{cccccc}
          1 & 0 & -\sqrt{2}+1 & \sqrt{2}-2 & \sqrt{2}+1 & \sqrt{2}+2  \\
          0 & 1 &     1       &      0     &     1      &     0
        \end{array} \right ] , \,\,
  J = \left [\begin{array}{cccccc}
     0 & 0 & 0 & 0 & 0  &  0 \\
     0 & 0 & 0 & 0 & 0  &  0 \\
     0 & 0 & 1 & 1 & 0  &  0 \\
     0 & 0 & 0 & 1 & 0  &  0 \\
     0 & 0 & 0 & 0 & -1 &  1 \\
     0 & 0 & 0 & 0 & 0  & -1
    \end{array} \right ]
\]
\[
 \mbox{and} \;\;\; Y^T =\, \frac{1}{4} \left [ \begin{array}{cccccc}
          0 & -4 & \sqrt{2}+2  & -\sqrt{2}-1&-\sqrt{2}+2 &-\sqrt{2}+1  \\
          4 &  0 &     0       &      1     &     0      &     -1
        \end{array} \right ] .
\]
The associated condition number of the eigenproblem of $P(\la)$ is
$\,k(P)=6.4183$. For $\eps = 0.3$ and $\ww = \{w_0, w_1, w_2,
w_3\} = \{ 0.1, 1, 1, 0 \}$, the matrix polynomial
\[
 Q(\la) \,=\,  I \la^3
             + \left [ \begin{array}{cc}
                              \ii \,0.3 & \sqrt{2}  \\
                             \sqrt{2}   & -\ii\, 0.3
                                  \end{array} \right ] \la^2
               + \left [ \begin{array}{cc}
                                    0 & -0.7  \\
                                  0.7 &  0
                                   \end{array} \right ] \la
               + \left [ \begin{array}{cc}
                                 0.01 & 0  \\
                                  0   & 0.03
                                   \end{array} \right ]
\]
lies on the boundary of $\mathcal{B} \left ( P ,0.3, \ww \right )$
and has $\mu = 0.5691+\ii\,0.0043$ as an eigenvalue. Then
$\,\min\limits_{\la \in \si(P)} \left | \mu - \la \right | =
|0.5691+\ii\,0.0043 - 1| = 0.4309 ,\,$ the upper bound of Theorem
\ref{thm: thmElsner} is $\,0.8554,\,$ and the upper bound of
Theorem \ref{thm: thmChu} is $\,3.8240$.

It is clear that the Elsner-like upper bound is tighter than the
upper bound of Theorem \ref{thm: thmChu} when $\|P(\mu)\|$ is
sufficiently small; this is the case in the above example, where
$\|P(0.5691+\ii\,0.0043)\|=1.0562$. In particular, if we define
the quantity
\[
     \Omega(P,\eps,\mu) \,=\,  \left \{  \begin{array}{l}
   | \det A_m | \ ( p \ k(P) )^{mn} \, (\eps\,w(|\mu|))^{mn-1},
   \; \mbox{when}\;\; p \ k(P) \, \eps \, w(|\mu|) \geq 1        \\
   | \det A_m | \ ( p \ k(P) )^{\frac{mn}{p}} \,
   ( \eps \, w(|\mu|) )^{\frac{mn}{p}-1} ,
   \; \mbox{when}\;\; p \ k(P) \, \eps \, w(|\mu|) < 1
 \end{array}   \right .   ,
\]
then it is straightforward to see that the bound of Theorem
\ref{thm: thmElsner} is better than the bound of Theorem \ref{thm:
thmChu} if and only if $\; \| P(\mu) \| \,<\,
\Omega(P,\eps,\mu)^{\frac{1}{mn-1}}$.

% -------------------------------------------------------

%%%%%%%%%%%%%%%%%%%%%%%%%%%%%%%%%%%%%%%%%%%%%%%%%
\end{document}